\documentclass[12pt]{article}

\usepackage{mathrsfs}
\usepackage{amsfonts}
\usepackage{tikz}
\usepackage{amssymb,bm}
\usepackage{amsmath}
\usepackage{amsthm}
\usepackage{bbding}
\usetikzlibrary{matrix,shapes,snakes}

\allowdisplaybreaks
\newtheorem{thm}{Theorem}[section]

\newtheorem{lem}[thm]{Lemma}
\newtheorem{cor}[thm]{Corollary}

\newtheorem{rmk}[thm]{Remark}

\def\proof{\vskip 1mm\noindent{\it Proof.}\quad}
\newcommand{\qedd}{\hspace*{\fill}$\Box$\medskip}

\renewcommand{\vec}[1]{\bm{#1}}

\def\Tr{\hbox{\rm{Tr}}}

\begin{document}

\title{The compositional inverse of a class of bilinear permutation
polynomials over finite fields of characteristic 2\thanks{Partially
supported by National Basic Research Program of China
(2011CB302400).}}

\author{Baofeng Wu\thanks{Key Laboratory of Mathematics
Mechanization, AMSS, Chinese Academy of Sciences,
 Beijing 100190,  China. Email: wubaofeng@amss.ac.cn}, Zhuojun Liu\thanks{Key Laboratory of Mathematics
Mechanization, AMSS, Chinese Academy of Sciences,
 Beijing 100190,  China. Email: zliu@mmrc.iss.ac.cn}}
 \date{}

\maketitle

\begin{abstract}
A class of bilinear permutation polynomials over a finite field of
characteristic 2 was constructed in a recursive manner recently
which involved some other constructions as special cases. We
determine the compositional inverses of them based on a direct sum
decomposition of the finite field. The result generalizes that in
[R.S. Coulter, M. Henderson, The compositional inverse of a class of
permutation polynomials over a finite field, Bull. Austral. Math.
Soc. 65 (2002) 521-526].\vskip .5em


\end{abstract}


\section{Introduction}\label{secintro}

Let $\mathbb{F}_{q}$ be the finite field with $q$ elements where $q$
is a prime or a prime power, and $\mathbb{F}_{q}[x]$ be the ring of
polynomials over $\mathbb{F}_{q}$. For any
$f(x)\in\mathbb{F}_{q}[x]$, it can induce a map from
$\mathbb{F}_{q}$ to itself. $f(x)$ is called a permutation
polynomial if the map induced by it is bijective. In fact, we need
only to consider polynomials of degree less than $q$ when talking
about permutation behavior of them. Clearly, under the operation of
composition of polynomials and subsequent reduction modulo
$(x^q-x)$, the set of all permutation polynomials over
$\mathbb{F}_{q}$ forms a group which is isomorphic to $\mathcal
{S}_q$, the symmetric group on $q$ letters. Hence for any
permutation polynomial $f(x)\in\mathbb{F}_{q}[x]$, there exists a
unique polynomial $f^{-1}(x)\in\mathbb{F}_{q}[x]$ such that
$f(f^{-1}(x))\equiv f^{-1}(f(x))\equiv x \mod (x^q-x)$. $f^{-1}$ is
called the compositional inverse of $f$ (or vice versa).

Discovering new classes of permutation polynomials is an old and
important problem due to their applicable value in cryptography,
coding theory and combinatorics. However, it is far from easy to do
this. There are only a few classes of permutation polynomials known.
See \cite{lidl2, lidl3} for a survey of this topic and \cite{wang1,
akbary1, charpin, xdhou, zbzha}, for example, for some recent
progresses.

Given a class of permutation polynomials, it seems to be an even
more difficult problem to find the class of permutation polynomials
that represent their compositional inverses. It was noted in
\cite{coulter} that only for permutation linear polynomials,
monomials and Dickson polynomials, the compositional inverses could
be explicitly determined. To the knowledge of the authors, this list
was enlarged in recent years and compositional inverses of the
following several other classes of permutation polynomials were
determined:

(1) Permutation polynomials of the form $x^rf(x^s)$ over
$\mathbb{F}_q$ where $s|(q-1)$. Permutation behavior of such
polynomials was studied in \cite{dwan, zieve, akbary2}, and their
compositional inverses were obtained in \cite{wang2}.

(2) The linearized polynomials over $\mathbb{F}_{q^n}$. A polynomial
over $\mathbb{F}_{q^n}$ of the shape
$L(x)=\sum_{i=0}^{n-1}a_ix^{q^i}$ is called a linearized polynomial.
It is well known that $L(x)$ is a permutation polynomial if and only
if the matrix
\[D_L=\begin{pmatrix}
a_0&a_1&\dots&a_{n-1}\\
a_{n-1}^q&a_0^q&\dots&a_{n-2}^q\\
\vdots&\vdots&&\vdots\\
a_1^{q^{n-1}}&a_2^{q^{n-1}}&\dots&a_0^{q^{n-1}}
\end{pmatrix}\]
is non-singular \cite{lidl1}. In \cite{wu} the authors found that
the compositional inverse of $L(x)$ can be represented by cofactors
of elements in the first column of $D_L$ (see \cite[Theorem
4.5]{wu}).

(3) The bilinear polynomial $x\left(\Tr(x)+ax\right)$ over
$\mathbb{F}_{q^n}$ where $q$ is even and $n$ is odd, proposed in
\cite{blok}. Its compositional inverse was determined in
\cite{coulter} which is of a complicated form.

In this paper, we focus on extending the list above. More
definitely, we want to replace (3) in the above list with a more
general case. A polynomial over $\mathbb{F}_{q^n}$ is called a DO
polynomial if it is of the shape
\[\sum_{0\leq i,\,j\leq n-1}a_{ij}x^{q^i+q^j}.\]
A polynomial of the shape $L_1(x)L_2(x)$ for two linearized
polynomials over $\mathbb{F}_{q^n}$ is obviously a DO polynomial,
which is called a bilinear polynomial in \cite{blok}. It was also
raised as a problem finding bilinear permutation polynomials in
\cite{blok}, which could be reduced to the problem of finding
bilinear permutation polynomials of the shape $xL(x)$ for a
linearized polynomial $L(x)$. The special class in (3) above was
constructed (see \cite[Theorem 5]{blok}) and the compositional
inverse class was obtained afterwards (see \cite[Theorem
1]{coulter}).

We notice that the class in (3) was generalized in \cite{chapuy} in
a recursive manner recently, thus it is a natural question how to
generalize its compositional inverse to the compositional inverse of
the generalized class of bilinear permutation polynomials. This is
not direct since the method in \cite{coulter} to obtain
compositional inverse is a ``guess and determine" one: verifying the
result after guessing it based on some experimental evidences.
However, after further studying properties of such permutation
polynomials, we can overcome this difficulty. The main idea of our
method is to decompose the finite field into a direct sum of two
subspaces, and represent the map induced by a univariate permutation
polynomial $f(x)$ by a bivariate permutation polynomial system $\vec
 f(y,z)=\left(f_1(y,z),f_2(y,z)\right)$. In the case $f(x)$ is a bilinear
 polynomial we consider, the corresponding bivariate polynomial
 system $\vec f(y,z)$ is of a triangular form: $f_1(y,z)$ is
 independent of the variable $z$. We can get the inverse polynomial
 system after overcoming the difficulty of determining the
 inverse of a permutation induced by a linearized polynomial on a component of the direct sum decomposition
 of the finite field. To summarize, we transform the problem of
 computing inverse of a non-linear map on the finite field to the
 problem of computing inverse of a linear map on a subspace of it,
 which seems much easier to solve.

 The rest of the paper is organized as follows. In Section 2 we
 recall some constructions of bilinear permutation polynomials and
 determine their compositional inverses. In Section 3 we explain our
 method to obtain the results. Concluding remarks are given in
 Section 4.


\section{Bilinear permutations and their compositional inverses}\label{secbilin}

We denote the trace map from $\mathbb{F}_{q^n}$ to $\mathbb{F}_{q}$
by $\Tr_{\mathbb{F}_{q^n}/\mathbb{F}_{q}}$ or $\Tr$ for simplicity
when it will not cause confusion, that is
\[\Tr(x)=\sum_{i=0}^{n-1}x^{q^i}, ~ x\in\mathbb{F}_{q^n}.\]
Throughout the rest of the paper we only consider finite fields of
characteristic 2.

We firstly recall the construction of a class of bilinear
permutation polynomials proposed in \cite{blok}.

\begin{thm}[\cite{blok}]\label{bipp}
Let $q$ be even and $n$ be odd. Then the polynomial
\[f(x)=x\left(\Tr(x)+ax\right)\]
is a permutation polynomial over $\mathbb{F}_{q^n}$ for all
$a\in\mathbb{F}_{q}\backslash\{0,1\}$.
\end{thm}

In \cite{chapuy} Laigle-Chapuy generalized this construction in a
recursive manner.

\begin{thm}[\cite{chapuy}]\label{gbipp}
Let $q$ be even and $n$ be odd. Assume $xL(x)$ is a bilinear
permutation polynomial over $\mathbb{F}_{q}$ for a linearized
polynomial $L(x)\in\mathbb{F}_{q}[x]$. Then the polynomial
\[F(x)=x\left(L(\Tr(x))+a\Tr(x)+ax\right)\]
is a bilinear permutation polynomial over $\mathbb{F}_{q^n}$ for any
$a\in\mathbb{F}_{q}^*$.
\end{thm}

Note that the polynomial in Theorem \ref{bipp} can be derived from
Theorem \ref{gbipp} by setting $L(x)=x$. In the following we will
propose the compositional inverse of $F(x)$ given in Theorem
\ref{gbipp}. Firstly we remark that in representing maps from
$\mathbb{F}_{q^n}$ to itself by polynomials over $\mathbb{F}_{q^n}$,
we sometimes distinguish $\frac{1}{x}$ with $x^{q^n-2}$. For
example, we use $\frac{1}{\Tr(x)}$ to represent
$\Tr(x)^{q^n-2}=\Tr(x)^{q-2}$. Besides, we sometimes use $x^{1/2}$
instead of $x^{q^n/2}$.

\begin{thm}\label{invgbipp}
Use the same notations as in Theorem \ref{gbipp} and let $q=2^m$ for
a positive integer $m$.  Assume the compositional inverse of $xL(x)$
is $g(x)\in\mathbb{F}_{q}[x]$. Then

\begin{eqnarray*}
   F^{-1}(x)
   &=&a^{2^{m-1}-1}x^{2^{nm-1}}+\left(g(\Tr(x))+a^{2^{m-1}-1}\sum_{k=1}^{\frac{n-1}{2}}x^{2^{(2k-1)m-1}}\right)\cdot\\
   &&\left(\frac{\Tr(x)}{g(\Tr(x))}+ag(\Tr(x))\right)^{q-1}\\[.1cm]
   &&+\sum_{j=0}^{m-2}a^{2^{j}-1}\left(\frac{\Tr(x)}{g(\Tr(x))}+ag(\Tr(x))\right)^{2^m-2^{j+1}}
   \left(\sum_{k=0}^{\frac{n-1}{2}}x^{q^{2k}}\right)^{2^j}.
\end{eqnarray*}

\end{thm}

\proof We proceed by directly verifying $F^{-1}(F(x))=x$ under
subsequent reduction modulo $(x^{q^n}-x)$. Note that for any
$x\in\mathbb{F}_{q^n}$,
\[F^{-1}(x)=\left(\frac{x}{a}\right)^{1/2}\]
when $\frac{\Tr(x)}{g(\Tr(x))}+ag(\Tr(x))=0$, and
\begin{eqnarray*}
  F^{-1}(x) &=&g\left(\Tr(x)\right)+\sum_{j=0}^{m-1}\frac{a^{2^j-1}}{\left(\frac{\Tr(x)}{g \left(\Tr(x)\right)}+ag\left(\Tr(x)\right)\right)^{2^{j+1}-1}}
\left(\sum_{k=0}^{\frac{n-1}{2}}x^{q^{2k}}\right)^{2^j}  \\
   &=&g\left(\Tr(x)\right)+\sum_{j=0}^{m-1}\frac{a^{2^j-1}}{\left(\frac{\Tr(x)}{g \left(\Tr(x)\right)}+ag\left(\Tr(x)\right)\right)^{2^{j+1}-1}}
\sum_{k=0}^{\frac{n-1}{2}}x^{2^{2km+j}}
\end{eqnarray*}
otherwise.

Firstly, it is obvious that
\begin{eqnarray*}
   \Tr\left(F(x)\right)&=& \Tr\left(xL(\Tr(x))+ax\Tr(x)+ax^2\right) \\
   &=& \Tr(x)L(\Tr(x))+ a\Tr(x)^2+a\Tr(x^2)\\
   &=&\Tr(x)L(\Tr(x))
\end{eqnarray*}
as $\Tr(1)=1$ and $\Tr(x^2)=\Tr(x)^2$. Thus
$$g\left(\Tr(F(x))\right)=g\left(\Tr(x)L(\Tr(x))\right)=\Tr(x)$$ since
$g(xL(x))=x$. Besides, for any $x\in\mathbb{F}_{q^n}$, it is clear
that $L(\Tr(x))+a\Tr(x)=0$ if and only if
$\frac{\Tr(F(x))}{g(\Tr(F(x)))}+ag(\Tr(F(x)))=0$ since
\[\frac{\Tr(F(x))}{g(\Tr(F(x)))}+ag(\Tr(F(x)))=\frac{\Tr(x)L(\Tr(x))}{\Tr(x)}+a\Tr(x).\]
When  $L(\Tr(x))+a\Tr(x)=0$, we have $F(x)=ax^2$ and hence
\[F^{-1}(F(x))=\left(\frac{x}{a}\right)^{1/2}\circ(ax^2)=x;\]
When $L(\Tr(x))+a\Tr(x)\neq0$, we have

\begin{eqnarray*}
  F^{-1}(F(x)) &=& \Tr(x)+
  \sum_{j=0}^{m-1}\frac{a^{2^j-1}}{\left(\frac{\Tr(x)L\left(\Tr(x)\right)}{\Tr(x)}+a\Tr(x)\right)^{2^{j+1}-1}}\cdot\\
&&\qquad~~~~~\sum_{k=0}^{\frac{n-1}{2}}\left(xL(\Tr(x))+ax\Tr(x)+ax^2\right)^{2^{2km+j}}\\
   &=& \Tr(x)+ \sum_{j=0}^{m-1}\frac{a^{2^j-1}}{\left(L(\Tr(x))+a\Tr(x)\right)^{2^{j+1}-1}}\cdot\\
&&\qquad~~\Bigg[\left(L(\Tr(x))+a\Tr(x)\right)^{2^j} \sum_{k=0}^{\frac{n-1}{2}}x^{2^{2km+j}}+a^{2^j}\sum_{k=0}^{\frac{n-1}{2}}x^{2^{2km+j+1}}\Bigg] \\
   &=& \Tr(x)+\sum_{j=0}^{m-1}\frac{a^{2^j-1}}{\left(L(\Tr(x))+a\Tr(x)\right)^{2^{j}-1}}\sum_{k=0}^{\frac{n-1}{2}}x^{2^{2km+j}}\\
  &&\qquad~ +\sum_{j=0}^{m-1}\frac{a^{2^{j+1}-1}}{\left(L(\Tr(x))+a\Tr(x)\right)^{2^{j+1}-1}}\sum_{k=0}^{\frac{n-1}{2}}x^{2^{2km+j+1}}\\
   &=&\Tr(x)+\sum_{k=0}^{\frac{n-1}{2}}x^{2^{2km}}+\sum_{k=0}^{\frac{n-1}{2}}x^{2^{(2k+1)m}}\\
   &=& \Tr(x)+\sum_{k=0}^{n-1}x^{q^k}+x^{q^n} \\
   &=&x.
\end{eqnarray*}
To summarize, we have $F^{-1}(F(x))=x$ for any
$x\in\mathbb{F}_{q^n}$.\qedd

\begin{cor}\label{invbipp}
Let $q=2^m$ and $n$ be odd. Let $f(x)$ be the bilinear permutation
polynomial defined in Theorem \ref{bipp}. Then
\begin{eqnarray*}
f^{-1}(x)&=&a^{2^{m-1}-1}x^{2^{nm-1}}+(1+a)^{2^{m-1}-1}\Tr(x)^{2^{m-1}}\\
&&+a^{2^{m-1}-1}\Tr(x)^{2^{m-1}(2^m-1)}\sum_{k=1}^{\frac{n-1}{2}}x^{2^{(2k-1)m-1}}\\
&&+\sum_{j=0}^{m-2}a^{2^j-1}(1+a)^{2^{m-1}+2^j-1}{\Tr(x)^{2^{m-1}-2^j}}
\left(\sum_{k=0}^{\frac{n-1}{2}}x^{q^{2k}}\right)^{2^j}.
\end{eqnarray*}

\end{cor}
\proof We let $L(x)=x$ and replace $a$ with $\frac{a}{1+a}$ in
Theorem \ref{gbipp}. Then
\begin{eqnarray*}
  F(x) &=& x\left(\frac{1}{1+a}\Tr(x)+\frac{a}{1+a}x\right) \\
   &=& \frac{1}{1+a}x\left(\Tr(x)+ax\right) \\
   &=&\frac{1}{1+a}f(x).
\end{eqnarray*}
Hence $f^{-1}(x)=F^{-1}\left(\frac{x}{1+a}\right)$. Then the result
is obtained from Theorem \ref{invgbipp} noting that
$(xL(x))^{-1}=x^{1/2}$.\qedd

\begin{rmk}
\cite[Theorem 1]{coulter} can be got from Corollary \ref{invbipp} by
replacing $a$ by $1+\alpha$ for
$\alpha\in\mathbb{F}_{q}\backslash\{0,1\}$. Though the compositional
inverse deduced from Corollary \ref{invbipp} is not totally the same
with that in \cite[Theorem 1]{coulter} in representation, they are
indeed the same polynomial after collections of terms.
\end{rmk}

Recently, Dempwolff and M\"uller constructed a new class of bilinear
permutation polynomials in \cite{demp} using trace maps over a tower
of finite fields in constructing translation planes of even order.

\begin{thm}[\cite{demp}]\label{bipptow}
Let $d_1,~d_2,~\ldots,~d_h,~n$ be all positive integers satisfying
that $d_1|d_2|\cdots|d_h|n$ and $\frac{n}{d_1}$ is odd. Let
$c_i\in\mathbb{F}_{2^{d_i}}^*$, $1\leq i\leq h$, such that
$\sum_{j=1}^ic_j\neq 0$ for all $i$. Choose $1\leq l<d_1$ with
$\gcd(2^{d_1}-1, 2^l+1)=1$ and $c_0\in\mathbb{F}_{2^{d_1}}^*$. Set
\[L_{h+1}(x)=\left(\sum_{i=1}^hc_i\right)x+\sum_{i=1}^hc_iT_{n:d_i}(x)+c_0T_{n:d_1}(x)^{2^l},\]
where $T_{n:d_i}$ denotes the trace map from $\mathbb{F}_{2^n}$ to
$\mathbb{F}_{2^{d_i}}$, $1\leq i\leq h$. Then
$F_{h+1}(x)=xL_{h+1}(x)$ is a bilinear permutation polynomial over
$\mathbb{F}_{2^n}$.
\end{thm}

Actually, we find that the polynomial $F_{h+1}(x)$ in Theorem
\ref{bipptow} can be obtained from Theorem \ref{gbipp} recursively
(though maybe it was discovered independently) in the following
manner. Set
\[L_1(x)=c_0x^{2^l}\in\mathbb{F}_{2^{d_1}}[x]\]
and
\[L_i(x)=L_{i-1}(T_{d_i:d_{i-1}}(x))+\left(\sum_{j=1}^{i-1}c_j\right)T_{d_i:d_{i-1}}(x)+\left(\sum_{j=1}^{i-1}c_j\right)x\in\mathbb{F}_{2^{d_i}}[x]\]
for $2\leq i\leq h$. Finally we set
\[L_{h+1}(x)=L_{h}(T_{n:d_{h}}(x))+\left(\sum_{j=1}^{h}c_j\right)T_{n:d_{h}}(x)+\left(\sum_{j=1}^{h}c_j\right)x\in\mathbb{F}_{2^{n}}[x]\]
From the transitivity  of the trace map, it is easy to verify that
$L_{h+1}(x)$ is just the one  defined in Theorem \ref{bipptow}. By
Theorem \ref{gbipp}, we directly obtain that $F_i(x)=xL_i(x)$ is a
permutation polynomial over $\mathbb{F}_{2^{d_i}}$ for $2\leq i\leq
h$ and $F_{h+1}(x)$ is a permutation polynomial over
$\mathbb{F}_{2^n}$, since $F_1(x)=xL_1(x)$ is a permutation
polynomial over $\mathbb{F}_{2^{d_1}}$ due to
$\gcd(2^{d_1}-1,~2^l+1)=1$. Hence $F^{-1}_{h+1}(x)$ can be derived
from Theorem \ref{invgbipp} inductively.

\begin{cor}\label{invbipptow}
Use notations the same as in Theorem \ref{bipptow} and the
discussions after it,  and let $d_{h+1}=n$. Assume $u(2^l+1)\equiv 1
\mod (2^{d_1}-1)$ and $a_i=\sum_{j=1}^{i-1}c_j$, $2\leq i\leq h+1$.
Then
\[F_1^{-1}(x)=\left(\frac{x}{c_0}\right)^u,\]
and for $2\leq i\leq h+1$,
\begin{eqnarray*}
   &&F_i^{-1}(x)= a_i^{2^{d_{i-1}-1}-1}x^{2^{d_i-1}}\\
   &&+\left(F_{i-1}^{-1}(T_{d_i:d_{i-1}}(x))+a_i^{2^{d_{i-1}-1}-1}\sum_{k=1}^{\frac{d_i/d_{i-1}-1}{2}}x^{2^{(2k-1)d_{i-1}-1}}\right)\cdot\\
&&\quad~\left(\frac{T_{d_i:d_{i-1}}(x)}{F_{i-1}^{-1}(T_{d_i:d_{i-1}}(x))}+a_iT_{d_i:d_{i-1}}(x)\right)^{2^{d_{i-1}}-1}\\[.1cm]
&&+\sum_{j=0}^{d_{i-1}-2}a_i^{2^j-1}
\left(\frac{T_{d_i:d_{i-1}}(x)}{F_{i-1}^{-1}(T_{d_i:d_{i-1}}(x))}+a_iT_{d_i:d_{i-1}}(x)\right)^{2^{d_{i-1}}-2^{j+1}}
\sum_{k=0}^{\frac{d_i/d_{i-1}-1}{2}}x^{2^{2kd_{i-1}+j}}.
\end{eqnarray*}
\end{cor}


\section{The method to obtain Theorem \ref{invgbipp}}\label{seccomp}

The method we get Theorem \ref{invgbipp} is not a ``guess and
determine" one, as a matter of fact. In this section we describe it
in detail.

Let $q$ be even and $n$ be odd. The for any $c\in\mathbb{F}_{q}$,
$\Tr(c)=c$, which implies that $\mathbb{F}_{q}\cap\ker \Tr=\{0\}$,
where $\ker \Tr$ is the kernel of the trace map from
$\mathbb{F}_{q^n}$ to $\mathbb{F}_{q}$. Since $\mathbb{F}_{q}$ and
$\ker\Tr$ are 1-dimensional and $(n-1)$-dimensional vector spaces
over $\mathbb{F}_{q}$, respectively, the following lemma is
straightforward.

\begin{lem}\label{decomp}
Let $q$ be even and $n$ be odd. Then
\[\mathbb{F}_{q^n}\cong\mathbb{F}_{q}\oplus\ker\Tr.\]
\end{lem}

\begin{rmk}
The map to establish the isomorphism in Lemma \ref{decomp} is
\begin{eqnarray*}
   \phi:~\mathbb{F}_{q^n}&\longrightarrow& \mathbb{F}_{q}\oplus\ker\Tr \\
   x&\longmapsto&(\Tr(x),x+\Tr(x)),
\end{eqnarray*}
and for $(y,z)\in\mathbb{F}_{q}\oplus\ker\Tr$,
$\phi^{-1}((y,z))=y+z$.
\end{rmk}

Now we consider the graph
\begin{center}
\begin{tikzpicture}
\matrix (m) [matrix of math nodes,row sep=4em,column sep=5em,minimum
width=2em]
  {
     \mathbb{F}_{q^n} & \mathbb{F}_{q^n} \\
    \mathbb{F}_{q}\oplus\ker\Tr & \mathbb{F}_{q}\oplus\ker\Tr \\};
  \path[-stealth]
    (m-1-1) edge node [left] {$\phi$} (m-2-1)
            edge node [above] {$F(x)$} (m-1-2)
    (m-2-1) edge node [above] {$\vec F(y,z)$} (m-2-2)
    (m-1-2) edge node [right] {$\phi$} (m-2-2);
\end{tikzpicture}
\end{center}
where $F(x)$ is the bilinear permutation polynomial defined in
Theorem \ref{gbipp}, and $\vec F(y,z)$ is a bivariate polynomial
system the map induced by which can make the graph commutative. Let
$y=\Tr(x)$ and $z=x+\Tr(x)$. Since
$$\Tr(F(x))=\Tr(x)L(\Tr(x))=yL(y)$$
and
\begin{eqnarray*}
  F(x)+\Tr(F(x)) &=&(y+z)L(y)+a(y+z)y+a(y+z)^2+yL(y)\\
   &=&az^2+(L(y)+ay)z,
\end{eqnarray*}
we have
\[\vec F(y,z)=\left(yL(y),~ az^2+(L(y)+ay)z\right).\]
$yL(y)$ is a permutation polynomial over $\mathbb{F}_{q}$ as stated,
thus $az^2+(L(y)+ay)z$ can induce a permutation of $\ker \Tr$ for
any $y\in\mathbb{F}_{q}$ since $\vec F(y,z)$ can induce a
permutation of $\mathbb{F}_{q}\oplus\ker\Tr$. In fact, $\vec F(y,z)$
is a bivariate polynomial system of the so-called triangular form,
so the inverse polynomial system can be obtained in case we can get
the inverse of the permutation of $\ker\Tr$.

\begin{lem}\label{invkerpp}
Let $q=2^m$ and $n$ be odd. Let $P_c(x)=x^2+cx$ for any
$c\in\mathbb{F}_q^*$. Then $P_c(x)$ can induce a permutation of
$\ker \Tr$ and the polynomial that can induce its inverse map is
\[P_c^{-1}(x)=\sum_{j=0}^{m-1}c^{-(2^{j+1}-1)}\left(\sum_{k=0}^{\frac{n-1}{2}}x^{q^{2k}}\right)^{2^j}.\]
(We call $P_c(x)$ a permutation polynomial over $\ker\Tr$ and
$P_c^{-1}(x)$ its compositional inverse.)
\end{lem}
\proof Obviously, for any $x\in\ker\Tr$, $P_c(x)\in\ker\Tr$.
Furthermore, $P_c(x)$ induces an $\mathbb{F}_2$-linear
transformation of the vector space $\ker\Tr$, the kernel of which is
$\{0\}$ since $P_c(x)=0$ implies $x=0$ or $x=c$ but
$c\not\in\ker\Tr$. Hence the linear transformation induced by
$P_c(x)$ is invertible. For any $x\in\ker\Tr$,
\begin{eqnarray*}
  P_c^{-1}(P_c(x)) &=&\sum_{j=0}^{m-1}c^{-(2^{j+1}-1)}\sum_{k=0}^{\frac{n-1}{2}}(x^2+cx)^{2^{2km+j}}  \\
   &=&\sum_{j=0}^{m-1}c^{-(2^{j+1}-1)}\sum_{k=0}^{\frac{n-1}{2}}x^{2^{2km+j+1}}+
   \sum_{j=0}^{m-1}c^{-(2^{j}-1)}\sum_{k=0}^{\frac{n-1}{2}}x^{2^{2km+j}}  \\
&=&\sum_{k=0}^{\frac{n-1}{2}}x^{2^{(2k+1)m}} +\sum_{k=0}^{\frac{n-1}{2}}x^{2^{2km}} \\
   &=& \Tr(x)+x \\
   &=&x.
\end{eqnarray*}
Hence $P_c^{-1}(x)$ is just the compositional inverse of $P_c(x)$
over $\ker\Tr$.\qedd

\begin{rmk}
In fact, $P_c^{-1}(x)$ in Lemma \ref{invkerpp} is not got fully by
guessing and determining. Since $P_c(x)$ is a linearized permutation
polynomial over $\ker\Tr$, we can assume
$P_c^{-1}(x)=\sum_{i=0}^{mn-1}d_ix^{2^i}$. However, we cannot expect
to get $P_c^{-1}(P_c(x))=x$ over $\mathbb{F}_{q^{n}}$ since $P_c(x)$
is not a permutation polynomial over $\mathbb{F}_{q^{n}}$.  Hence we
expect to get
$P_c^{-1}(P_c(x))=x+\Tr(x)=x^{2^m}+x^{2^{2m}}+\cdots+x^{2^{(n-1)m}}$
as we will be limited to $\ker\Tr$, which leads to a system of
equations
\[\left\{\begin{aligned}
&d^2_{i-1}+cd_i=0&\text{if}~i\neq km\\
& d^2_{i-1}+cd_i=1&\text{if}~i= km
\end{aligned}\right.,~0\leq i\leq mn-1,~1\leq k\leq n-1.\]
Solving this system we obtain $d_i=c^{-(2^{j+1}-1)}$ when $k$ is
even and $d_i=0$ when $k$ is odd for $i=km+j$, $0\leq j\leq m-1$,
$0\leq k\leq n-1$.
\end{rmk}

Let $\vec F^{-1}(y,z)$ be the bivariate polynomial system such that
for any $(y,z)\in\mathbb{F}_{q}\oplus\ker\Tr$, $\vec F^{-1}\circ\vec
F(y,z)=(y,z)$. Now for $(Y,Z)\in\mathbb{F}_{q}\oplus\ker\Tr$, we
assume
\[\left\{\begin{aligned}
&yL(y)=Y\\
&az^2+(L(y)+ay)z=Z .
\end{aligned}\right.\]
Let $g(x)\in\mathbb{F}_{q}[x]$ be the compositional inverse of the
permutation polynomial $xL(x)$ over $\mathbb{F}_{q}$. If
$L(g(Y))+ag(Y)=0$, i.e. $\frac{Y}{g(Y)}+ag(Y)=0$ (since
$g(Y)L(g(Y))=Y$), we get
\[\left\{\begin{aligned}
y&=g(Y)\\
z&=\left(\frac{Z}{a}\right)^{1/2} .\end{aligned}\right.\] This is
equivalent to say
\[\vec F^{-1}(Y,Z)=\left(g(Y),~\left(\frac{Z}{a}\right)^{1/2}\right)=\left(\left(\frac{Y}{a}\right)^{1/2},~\left(\frac{Z}{a}\right)^{1/2}\right);\]
If $L(g(Y))+ag(Y)\neq0$, we get
\[\left\{\begin{aligned}
y&=g(Y)\\
z&=P_{L(y)/a+y}^{-1}\left(\frac{Z}{a}\right)=P_{Y/(ag(Y))+g(Y)}^{-1}\left(\frac{Z}{a}\right)
,\end{aligned}\right.\] where $P_c(x)$ is the polynomial defined in
Lemma \ref{invkerpp}. By Lemma \ref{invkerpp} we have
\begin{eqnarray*}
   z&=&\sum_{j=0}^{m-1}\left(\frac{Y}{ag(Y)}+g(Y)\right)^{-(2^{j+1}-1)}\sum_{k=0}^{\frac{n-1}{2}}\left(\frac{Z}{a}\right)^{2^{2km+j}}\\
   &=&\sum_{j=0}^{m-1}\frac{a^{2^{j}-1}}{\left(\frac{Y}{g(Y)}+ag(Y)\right)^{2^{j+1}-1}}\sum_{k=0}^{\frac{n-1}{2}}Z^{2^{2km+j}}.
\end{eqnarray*}
Hence
\[\vec
F^{-1}(Y,Z)=\left(g(Y),~\sum_{j=0}^{m-1}\frac{a^{2^{j}-1}}{\left(\frac{Y}{g(Y)}+ag(Y)\right)^{2^{j+1}-1}}\sum_{k=0}^{\frac{n-1}{2}}Z^{2^{2km+j}}\right).\]
 Now we go back to $\mathbb{F}_{q^n}$ by
setting $X=Y+Z$ where $Y=\Tr(X)$ and $Z=X+\Tr(X)$. We get when
$\frac{\Tr(X)}{g(\Tr(X))}+ag(\Tr(X))=0$,
\[F^{-1}(X)=\left(\frac{\Tr(X)}{a}\right)^{1/2}+\left(\frac{X+\Tr(X)}{a}\right)^{1/2}=\left(\frac{X}{a}\right)^{1/2},\]
and when $\frac{\Tr(X)}{g(\Tr(X))}+ag(\Tr(X))\neq0$,
\begin{eqnarray*}
   F^{-1}(X)&=&g(\Tr(X))+\sum_{j=0}^{m-1}\frac{a^{2^{j}-1}}{\left(\frac{\Tr(X)}{g(\Tr(X))}+ag(\Tr(X))\right)^{2^{j+1}-1}}\sum_{k=0}^{\frac{n-1}{2}}(X+\Tr(X))^{2^{2km+j}}  \\
   &=&g(\Tr(X))+\sum_{j=0}^{m-1}\frac{a^{2^{j}-1}}{\left(\frac{\Tr(X)}{g(\Tr(X))}+ag(\Tr(X))\right)^{2^{j+1}-1}}\sum_{k=0}^{\frac{n-1}{2}}X^{2^{2km+j}}\\[.1cm]
   &&\qquad\quad~~+\frac{n+1}{2}\sum_{j=0}^{m-1}\frac{a^{2^{j}-1}\Tr(X)^{2^j}}{\left(\frac{\Tr(X)}{g(\Tr(X))}+ag(\Tr(X))\right)^{2^{j+1}-1}}\\
&=&g(\Tr(X))+\sum_{j=0}^{m-1}\frac{a^{2^{j}-1}}{\left(\frac{\Tr(X)}{g(\Tr(X))}+ag(\Tr(X))\right)^{2^{j+1}-1}}\left(\sum_{k=0}^{\frac{n-1}{2}}X^{q^{2k}}\right)^{2^j}
\end{eqnarray*}
since
\begin{eqnarray*}
   &&\sum_{j=0}^{m-1}\frac{a^{2^{j}-1}\Tr(X)^{2^j}}{\left(\frac{\Tr(X)}{g(\Tr(X))}
   +ag(\Tr(X))\right)^{2^{j+1}-1}}\\[.2cm]
   &=&   \frac{\frac{\Tr(X)}{g(\Tr(X))}+
   ag(\Tr(X)}{a}\cdot\sum_{j=0}^{m-1}\frac{1}{\frac{\Tr(X)^{2^j}}{a^{2^j}g(\Tr(X))^{2^{j+1}}}+\frac{a^{2^j}g(\Tr(X))^{2^{j+1}}}{\Tr(X)^{2^j}}}
   \\[.2cm]
   &=&\frac{\frac{\Tr(X)}{g(\Tr(X))}+
   ag(\Tr(X)}{a}\cdot\Tr_{\mathbb{F}_{2^m}/\mathbb{F}_{2}}\left(\frac{1}{\frac{\Tr(X)}{ag(\Tr(X))^2}+\frac{ag(\Tr(X))^2}{\Tr(X)}}\right)\\
   &=&0
\end{eqnarray*}
using the following lemma.

\begin{lem}
Let $r$ be a positive integer. Then for any $e\in\mathbb{F}_{2^r}$,
$$\Tr_{\mathbb{F}_{2^r}/\mathbb{F}_{2}}\left(\frac{1}{e+e^{-1}}\right)=0.$$
$($Recall that we distinguish $\frac{1}{0}$ and $0^{-1}$ with
$0^{2^r-2}=0$.$)$
\end{lem}
\proof Since
\[\frac{1}{e+e^{-1}}=\frac{e}{e^2+1}=\frac{e+1+1}{(e+1)^2}=\frac{1}{e+1}+\frac{1}{(e+1)^2},\]
the result is straightforward to get.\qedd

Finally $F^{-1}(X)$ can be derived by interpolation, that is
\begin{eqnarray*}
   F^{-1}(X)&=& \left(\frac{X}{a}\right)^{1/2}\left[1+\left(\frac{\Tr(X)}{g(\Tr(X))}+ag(\Tr(X))\right)^{q-1}\right] \\
   &&+\left[g(\Tr(X))+\sum_{j=0}^{m-1}\frac{a^{2^{j}-1}}{\left(\frac{\Tr(X)}{g(\Tr(X))}+ag(\Tr(X))\right)^{2^{j+1}-1}}
   \left(\sum_{k=0}^{\frac{n-1}{2}}X^{q^{2k}}\right)^{2^j}\right]\cdot\\
   &&\left(\frac{\Tr(X)}{g(\Tr(X))}+ag(\Tr(X))\right)^{q-1}\\
   &=&\left(\frac{X}{a}\right)^{1/2}+\left[\left(\frac{X}{a}\right)^{1/2}+g(\Tr(X))+a^{2^{m-1}-1}\left(\sum_{k=0}^{\frac{n-1}{2}}X^{q^{2k}}\right)^{2^{m-1}}\right]\cdot\\
   &&\left(\frac{\Tr(X)}{g(\Tr(X))}+ag(\Tr(X))\right)^{q-1}\\[.1cm]
   &&+\sum_{j=0}^{m-2}a^{2^{j}-1}\left(\frac{\Tr(X)}{g(\Tr(X))}+ag(\Tr(X))\right)^{2^m-2^{j+1}}
   \left(\sum_{k=0}^{\frac{n-1}{2}}X^{q^{2k}}\right)^{2^j}\\
   &=&a^{2^{m-1}-1}X^{1/2}+\left(g(\Tr(X))+a^{2^{m-1}-1}\sum_{k=1}^{\frac{n-1}{2}}X^{2^{(2k-1)m-1}}\right)\cdot\\
   &&\left(\frac{\Tr(X)}{g(\Tr(X))}+ag(\Tr(X))\right)^{q-1}\\[.1cm]
   &&+\sum_{j=0}^{m-2}a^{2^{j}-1}\left(\frac{\Tr(X)}{g(\Tr(X))}+ag(\Tr(X))\right)^{2^m-2^{j+1}}
   \left(\sum_{k=0}^{\frac{n-1}{2}}X^{q^{2k}}\right)^{2^j}.
\end{eqnarray*}
Hence Theorem \ref{invgbipp} is obtained.


\section{Concluding remarks}\label{conclu}

In this paper, we derive the compositional inverse of a class of
bilinear permutation polynomials over a finite field of
characteristic 2, which can serve as a new class of permutation
polynomials over finite fields. The main observation we make is the
special structure of the maps induced by this class of permutation
polynomials: they can be represented by a class of bivariate
polynomial systems which are of triangular shapes after decomposing
the finite field into a direct sum of a subfield and the kernel
space of the trace map. In fact, the class of bilinear polynomials
in Theorem \ref{bipptow} can also be represented by multivariate
polynomial systems of triangular shapes considering the
decomposition
\[\mathbb{F}_{2^n}\cong\mathbb{F}_{2^{d_1}}\oplus\ker T_{d_2:d_1}\oplus\ker T_{d_3:d_2}\oplus\cdots\oplus\ker T_{n:d_h}.\]
Indeed, the corresponding polynomial systems are of the form
\[\left(~\begin{array}{l}
x_0L_1(x_0),\\[.1cm]
x_1[L_1(x_0)+c_1x_0]+c_1x_1^2, \\[.1cm]
x_2[L_2(x_0+x_1)+(c_1+c_2)(x_0+x_1)]+(c_1+c_2)x_2^2,\\[.1cm]
\quad\vdots\\[.1cm]
x_i\left[L_i\left(\sum_{j=0}^{i-1}x_j\right)+\left(\sum_{j=1}^{i}c_j\right)\left(\sum_{j=0}^{i-1}x_j\right)\right]+
\left(\sum_{j=1}^{i}c_j\right)x_i^2,\\[.1cm]
\quad\vdots\\[.1cm]
x_h\left[L_h\left(\sum_{j=0}^{h-1}x_j\right)+\left(\sum_{j=1}^{h}c_j\right)\left(\sum_{j=0}^{h-1}x_j\right)\right]+
\left(\sum_{j=1}^{h}c_j\right)x_h^2
\end{array}~\right).\]
Inverse polynomial systems of them can also be computed inductively
using Lemma \ref{invkerpp}, which will lead to Corollary
\ref{invbipptow} after lifted back to univariate polynomials.


\end{document}